\documentclass[9pt]{amsart}
\usepackage{amssymb}
\theoremstyle{plain}

\theoremstyle{definition}
\newtheorem{definition}{Definition}
\theoremstyle{remark}
\newtheorem{remark}{Remark}
\numberwithin{equation}{section}
\begin{document}
\title{Analytic hypoellipticity for $\square_b + c$ on the Heisenberg group: an
$L^2$ approach}
\author{David S. Tartakoff}
\address{Department of Mathematics, Statistics and
Computer Science\\
					University of Illinois at Chicago \\ 851 W. Morgan St. \\ Chicago
IL  60607 USA}
\email{dst@uic.edu}
%
%
\begin{abstract}
In an interesting note, E.M. Stein observed some 20 years ago that
while the Kohn Laplacian $\square_b$ on functions is neither
locally solvable nor (analytic) hypoelliptic, the addition of a
non-zero complex constant reversed these conclusions at least on the
Heisenberg group, and Kwon reproved and generalized this result
using the method of concatenations. Recently Hanges and Cordaro
have studied this situation on the Heisenberg group in detail. Here we give a
purely
$L^2$ proof of Stein's result using the author's now
classical construction of $(T^p)_\phi = \phi T^p +\ldots,$ where
$T$ is the 'missing direction' on the Heisenberg group. 
\end{abstract}
\maketitle

\section{Introduction}

As had been known for many years, the Kohn Laplacian behaves
very differently on functions and on forms. This is clearly stated
when the underlying manifold is the (complex) Heisenberg group,
which we now define. 

On $R^{2n-1} \sim C^{n-1}\times R$ we consider the complex
vector fields 
$$L_j = {\partial\over \partial z_j} - {i\over
2}\overline{z_j}{\partial\over \partial t}$$
which have the commutation relations 
$$[L_j, L_k] = 0, \quad [{L_j}, \overline{L_k}] =
\delta_{jk}i{\partial\over
\partial t}= -\delta_{jk}T.$$

The Kohn Laplacian is defined (on functions) by the
expression 
$$\square_b = \sum_j
\overline{L_j}^*\,\overline{L_j}$$
and on forms by 
$$\square_b = 
\overline{\partial}_b^*{\overline\partial_b} +
{\overline\partial_b}\overline{\partial}_b^*$$
where $\overline{\partial}_b$ on is the usual extension to forms of
the operator 
$$v\rightarrow \overline\partial_b v = (\overline{L}_1v, \ldots
\overline{L}_{n-1}v)$$ to form a complex: $\overline{\partial}^2
= 0$, and 
${}^*$ denotes $L^2$ adjoint. (These formulas are valid quite
generally on any CR manifold $M$  where the $\{L_j\}$ form a
basis for $T^{1,0}$ in the splitting $CTM = T^{1,0} + T^{0,1} +
F,$  with $T^{0,1} = \overline T^{1,0}$ and
${\hbox{dim}_R}(F=\overline{F}) = 1.$) 

Thus on the Heisenberg group (i.e., with the vector fields $L_j$ as
above), on functions 
$$\square_b = -\sum_j {L_j}{\overline{L}_j}.$$

Now while for a strictly pseudo-convex CR manifold in general, and
on the Heisenberg group in particular, $\square_b$ is hypoelliptic
and even analytic hypoelliptic on forms of degree at least one, this is not true on
functions and this operator is not even locally solvable. 

However in a striking note in 1982 (\cite{Stein}), E.M. Stein
showed that the addition of an arbitrary non-zero complex
constant to
$\square_b$ reverses this conclusion and $\square_b + c$ becomes
locally solvable and (analytic-~)hypoelliptic. His proof uses
explicit kernels. 

Soon after, Kwon (\cite{Kwon}) published a paper which, using the method of
concatenations of F. Treves, extended the
result to a non-zero real analytic {\it function}
$c(z,\overline{z}, t )$ instead of a constant. Recently Cordaro and Hanges
have been studying this situation in detail (\cite{Cordaro-Hanges}).

Here we give a totally elementary proof of the result on the
Heisenberg group.  

\section{The {\it a priori} estimate}

The starting point of any local solvability or hypoellipticity proof is
a suitable {\it a priori} estimate. For the operator 
$$\square_b + c= -\sum_j {L_j}{\overline{L}_j} + c,$$
we have 
$$((\square_b+c)v,v)_{L^2} = \sum_j\|\overline{L_j} v\|_{L^2}^2 + c
\|v\|_{L^2}^2,$$
so that if 1) $\Re{c}>0$ we may assert
$$|\Re((\square_b+c)v,v)_{L^2}| \geq \sum_j\|\overline{L_j} v\|_{L^2}^2 +
|\Re{c}| \|v\|_{L^2}^2,$$
while we always have 
$$|\Im((\square_b+c)v,v)_{L^2}|= |\Im{c}|\|w\|_{L^2}^2$$
and
$$|\Re((\square_b+c)v,v)_{L^2}| \geq \sum_j\|\overline{L_j} v\|_{L^2}^2 -
|\Re{c}| \|v\|_{L^2}^2,$$
and so if $|\Re{c}|<\tilde{d}|\Im{c}|$ (and $\Re{c}<0$):
$$|\Im((\square_b+c)v,v)_{L^2}| +
\tilde{d}|\Re((\square_b+c)v,v)_{L^2}|
\geq \tilde{d}\sum_j
\|\overline{L_j}v\|_{L^2}^2+(1-\tilde{d})\|v\|_{L^2}^2.$$
Thus for all $c$ except $c<0,$ and suitable $d=d(c),$ 
$$|((\square_b+c)v,v)_{L^2}| \geq C^\prime_c(\sum_j
\|\overline{L_j}v\|_{L^2}^2+\|v\|_{L^2}^2)$$
uniformly in $v\in C_0^\infty.$

Finally, since the spectrum of $\square_b$ is discrete, the case
$c<0$ leads to a finite dimensional kernel of $\square_b+c$ which
may be handled by a norm of negative order. Thus for any
complex $c\neq 0$ there exists a constant $C_c$ such that for all
smooth
$v$ of compact support in a neighborhood of the origin, 
$$\sum_j \|\overline{L_j}v\|_{L^2}^2+\|v\|_{L^2}^2 \leq C_c 
\{|((\square_b+c)v,v)_{L^2}| + \|v\|_{-1}^2\}.$$
\begin{remark} The estimate is also valid, with essentially the same derivation,
if the non-zero constant $c$ added to $\square_b$ is replaced by any
elliptic pseudodifferential operator of arbitrary order $s<2$, with the
$L^2$ norm on the left replaced by the norm in $H^s$ and the second
norm on the right by any relatively compact norm, for instance, the
norm $\|\cdot\|_{H^{s-1}}.$
\end{remark}

\section{The definition and properties of $(T^p)_\phi$}

There are several ways to localize derivatives - when the derivatives
are measured by vector fields over which we have maximal control
in the estimates, or their conjugates, simple multiplication by a
cut-off function of H\"ormander-Ehrenpreis type will suffice: 
\begin{definition} There exists a constant $C=C(n)$ depending only
on the dimension $n$ with the following properties: given two
bounded open set $\Omega_j$ with $\overline\Omega_1\subset
\Omega_2,$ and separation $d=dist(\Omega_1,\Omega_2^c),$ and
any natural number $N,$ there exists $\Psi_N \in C_0^\infty
(\Omega_2),$ identically equal to one on $\overline{\Omega_1},$
and such that for $|\alpha|\leq 3N,$
$$|D^\alpha \Psi_N| \leq C({C\over d})^{|\alpha|}N^{|\alpha|},\quad
|\alpha|\leq 3N.$$
\end{definition}
The construction of such functions is just a convolution of
$N$ identical `bump' functions with the characteristic function of a
set `midway' between the two sets $\Omega_j.$  The factor $3$ is of
course arbitrary - any other fixed multiple (independent of
$N$) would work as well, with a different choice of $C(n).$

The rough idea, elaborated below, is that when instead of $v,$
the function $\Psi_N \overline{L}^\alpha u$ is subjected to the
{\it a priori} estimate, brackets of $\square_b$ with
$\Psi_N\overline{L}^\alpha$ will give two kinds of terms - but at
least those where the derivatives in the differential operator land on
$\Psi_N$ will introduce a factor of $CN/d$ but decrease $\alpha$
by one (counting the derivatives in $\square_b$ as $\overline{L}$'s
for the moment), and after $N$ iterations the factor will be
$(CN/d)^N$ which is bounded by Stirling's formula by
$(C^\prime)^NN!$ leading to analyticity.) 

However when the derivative whose high powers (here $T^p$)
are being localized is not one optimally handled by the {\it a
priori} estimate, this procedure is not sufficient. For incurring
a derivative on
$\Psi_N,$ and hence a factor of $CN/d,$ is not offset by `gaining'
one of the $\overline{L}$'s - since now we do not have a large
supply of $\overline{L}$'s to use up, and to convert a $T$ to {\it
two} $\overline{L}$'s (actually one $L$ and one $\overline{L}$)
is too costly - it would introduce {\it two} factors of $CN/d$ to
decrease $p$ by one, leading to $C^NN!^2,$ the second Gevrey class
($|D^ru| \leq C^{r+1}r!^2$) and not the analytic, or first, Gevrey
class. 

To localize $T^p,$ since 
$$[L_j,\Psi T^p] = (L_j \Psi)T^p,$$
on the Heisenberg group, we use the excellent commutation
relations, in particular that $[L_j, \overline{L_k}] 
= \delta_{jk}T,$ all other brackets being zero, 
to construct
$$T_\Psi = \Psi T + \sum(\overline{L_j}\Psi)L_j -
\sum({L_j}\Psi)\overline{L_j}$$
which evidently satisfies 
$$[L_k, T_\Psi] = \sum_j(L_k\overline{L_j}\Psi)L_j -
\sum_j(L_k{L_j}\Psi)\overline{L_j}$$
and
$$[\overline{L_k}, T_\Psi] =
\sum_j(\overline{L_k}\overline{L_j}\Psi)L_j -
\sum_j(\overline{L_k}{L_j}\Psi)\overline{L_j}.$$

At first sight this may not appear to be much of an improvement
over brackets with $\Psi T$ alone, but if we observe that here the
number of derivatives on the right on $\Psi$ should not be taken as
{\it double} the loss of $T$ derivatives (which would again lead
to the second Gevrey class) but rather as {\it one more} than the
loss of $T$ derivatives, we are encouraged to try to generalize this
construction for higher powers of $T.$ 

Order matters here, and after much trial and error it was found that
the following is one eminently satisfactory localization of $T^p:$ 
$$(T^p)_\Psi = \sum_{|\alpha + \beta|\leq p}
{(-1)^{|\alpha|}(L^\alpha \overline{L}^\beta
\Psi)\over \alpha!\beta!} L^\beta \overline{L}^\alpha T^{p-|\alpha +
\beta|}.$$
For this localization we have 
$$[L,(T^p)_\Psi]\equiv 0 $$
and
$$[\overline{L_k},(T^p)_\Psi]\equiv (T^{p-1})_{T\Psi}\circ
\overline{L_k}$$ modulo $C^p$ terms of the form 
$${(ML^\alpha \overline{L}^\beta
\Psi)\over \alpha!\beta!} L^\beta \overline{L}^\alpha $$
with $|\alpha +\beta|=p$
and $M=L$ or $M=\overline{L};$
we will write this error as 
$$C^{P+1}\Psi^{(p+1)}M^p/p!$$
with each instance of $M=L$ or $M=\overline{L}.$

If we apply these elegant relations to our solution via the {\it a
priori} estimate for $v=(T^p)_\Psi u$ and $(\square_b+c)u=f\in
C^\omega,$ we have
$$\sum_j \|\overline{L_j}(T^p)_\Psi
u\|_{L^2}^2+\|(T^p)_\Psi u\|_{L^2}^2 \leq  $$
$$\leq C\{|((\square_b+c)(T^p)_\Psi
u,(T^p)_\Psi u)_{L^2}| +
\|(T^p)_\Psi u\|_{-1}^2\}$$
$$\lesssim |((T^p)_\Psi f, (T^p)_\Psi u)_{L^2}|+
\|(T^p)_\Psi u\|_{-1}^2+$$
$$+ |(\sum_k
[L_k\overline{L_k},(T^p)_\Psi] u,(T^p)_\Psi u)_{L^2}|,$$
where the notation $A\lesssim B$ denotes $A \leq C B$ with a
constant $C$ depending only on the dimension $n.$

Thus, expanding the bracket, 

$$[L_k\overline{L_k},(T^p)_\Psi] u = 
L_k[\overline{L_k},(T^p)_\Psi] u +
[L_k,(T^p)_\Psi] \overline{L_k}u = 
L_k \circ (T^{p-1})_{T\Psi}\circ \overline{L_k} u$$
$$=(T^{p-1})_{T\Psi}\circ L_k \overline{L_k} u$$
modulo $C^p$ terms of the form $\Psi^{(p+1)}M^{p+1}/p!$
and in fact, due to the simple form of $\square_b,$ upon summing 
we recover $\square_b,$ though this will not help us.

Suffice it to say that from this expansion, upon iteration and with
the weighted Schwarz inequality we may write: 
$$\sum_j \|\overline{L_j}(T^p)_\Psi
u\|_{L^2}^2+\|(T^p)_\Psi u\|_{L^2}^2 \lesssim 
\sum_{0\leq q\leq p}C^q\|(T^{p-q})_{T^q\Psi}
f\|_{L^2}^2 +$$
$$+ \sum_{1\leq q\leq
p}C^q\|\overline{L_j}(T^{p-q})_{T^q\Psi} u\|_{L^2}^2
+\sum_{1\leq q\leq p}C^q\|(T^{p-q})_{T^q\Psi}
u\|_{L^2}^2 +$$
$$+\sum_{0\leq q\leq p}C^q\|(T^{p-q})_{T^q\Psi}
u\|_{-1}^2+\sum_{q\geq
0}\left(C^q\|\Psi^{(p+1)}M^{p-q+1}u\|/(p-q)!\right)^2$$

Now if the support of $\Psi$ is small enough, the $-1$ norm is less
than a small multiple of the $L^2$ norm, so the next to last term on
the right will be absorbed by the second term on the left (and the
third term on the right) - the trickier term is the last on the right. 

For this term we will introduce a {\it new} cut-off function,
$\tilde{\Psi},$ equal to one on the support of $\Psi,$ and with
essentially the same growth of its derivatives (they can be made the
same by taking $d$ half as large originally and nesting three open
sets instead of two). Then bring the high derivatives of $\Psi$ out of
the norm bounded by $(C/d)^{p+1}N^{p+1} \sim
(\tilde{C}/d)^{p+1}p!$ when $N$ is comparable to $p$. This
iteration shows that to bound high $T$ derivatives applied to $u$
locally by the corresponding factorial it will suffice to bound
derivatives of the same order but in the `directions' $L$ and
$\overline{L}$ by the same factorials. 

\section{Mixed powers of  $L_k$ and $\overline{L_j}$}

In this section we consider mixed $L$ and $\overline{L}$
derivatives of u. We claim that judicious integration by parts will
reduce us to consideration of pure powers of $\overline{L}$ modulo
{\it half as many} powers of $T.$ For this, we need to treat
expressions such as 
$$(L^\alpha \overline{L}^\beta w, L^\alpha \overline{L}^\beta
w)_{L^2} - (\overline{L}^\alpha \overline{L}^\beta w,
\overline{L}^\alpha
\overline{L}^\beta w)_{L^2}= \pm
 ([\overline{L}^\alpha,L^\alpha] \overline{L}^\beta w,
\overline{L}^\beta w)_{L^2}$$
for $w\in C_0^\infty.$

Now an elementary calculation expresses
$[\overline{L}^\alpha,L^\alpha]$ as a sum 
$$[\overline{L}^\alpha,L^\alpha]=\sum_{0\neq \alpha^\prime\leq
\alpha} {\alpha\choose\alpha^\prime}^2
\alpha^\prime!T^{|\alpha^\prime
|}L^{\alpha-\alpha^\prime}\overline{L}^{\alpha-\alpha^\prime},$$
and then we merely integrate $L^{\alpha-\alpha^\prime}$ and
approximately {\it half} of the $T$'s by parts in the inner
product. 

Howerver, to sum things up so far, we note (taking $f=0$ for
simplicity) that iterating the {\it a priori} estimate above and
bringing $\Psi^{(p+1)}$ out of the norm, we arrive, for $p$
comparable to $N,$ at 

$${\sum_j \|\overline{L_j}(T^p)_\Psi
u\|_{L^2}+\|(T^p)_\Psi u\|_{L^2}\over p!} \leq $$
$$\leq (\tilde{C}/d)^pp!\sup_{p\geq q+2r\geq
0}C^qr!{\|\tilde{\Psi}\overline{L}^{p-q-2r+1}T^ru\|\over
(p-q)!}$$ or, more suggestively, in the region where $\Psi
= \Psi_1 \equiv 1,$ and writing $d_1$ for that $d,$
$${\sum_j \|\overline{L_j}T^pu\|_{L^2(\Psi_1\equiv
1)}+\|T^pu\|_{L^2(\Psi_1\equiv 1)} \over p!}\leq $$
$$\leq (C/d_1)^p\sup_{p\geq q+2r\geq
0}C^q{\|\overline{L}{\Psi_2}\overline{L}^{p-q-2r}T^ru\|\over
(p-q-r)!}$$
with $\Psi_2\equiv 1$ on the support of $\Psi_1.$ Note that the
order of $T$ derivatives here can not exceed half of the order on the
left hand side, and the same will be true even after brackets when
brackets of $L$'s and $\overline{L}$'s produce additional $T$'s. 

\section{Brackets with $\overline{L}^\alpha$ (and a few
$T$'s and perhaps one $L$)}

The analysis of powers of $\overline{L}$ is relatively
straightforward, since we have a supply of `good' vector fields to
use over and over in the {\it a priori} estimate. All that happens is
that in the bracket with $\square_b (+c),$ when $L$ meets
$\overline{L}^{p-q-2r}$ there may be as many as ${p-q-2r}$
copies of $T$ and the total number of 
$\overline{L}$'s will be decreased by two: 
$$([L_k \overline{L_k}, \Psi_2 \overline{L}^s]T^m u, \Psi_2
\overline{L}^sT^m u)_{L^2}\sim (\Psi_2
\overline{L_k}\overline{L}^sT^m u, \Psi_2
^\prime\overline{L}^sT^m u)_{L^2}+$$
$$+(\Psi_2^\prime
\overline{L}^sT^m u, \overline{L_k}\Psi_2
\overline{L}^sT^m u)_{L^2}
+s(\Psi_2 \overline{L}^sT^{m+1} u,
\Psi_2
\overline{L}^sT^m u)_{L^2}
$$
(and some terms with two derivatives on $\Psi,$ and a drop
of two $\overline{L}$'s - in fact we may here safely ignore
derivatives on $\Psi$, and for instance we have interchanged
the two localizing functions in the first term on the right). 

To handle the last term we do not try to imagine half a power
of $T$ belonging to each side of the inner product, though we
could do so using pseudodifferential operators, but rather do
the paradoxical manoever: 
$$s(\Psi_2 \overline{L}^sT^{m+1} u, 
\Psi_2 \overline{L}^sT^m u)_{L^2} \sim  
s(\Psi_2 \overline{L}^{s-1}T^{m+1} u, 
{L}\Psi_2 \overline{L}^{s}T^m u)_{L^2}$$
$$\sim s\|\overline{L}\Psi_2 \overline{L}^{s-2}T^{m+1}
u\|_{L^2}\left(\| \overline{L}\Psi_2
L\overline{L}^{s-1}T^mu\|_{L^2}
+ \| \Psi_2\overline{L}^{s-1}T^{m+1})u\|_{L^2}\right)\sim$$
$$\sim C_\varepsilon (s\|\overline{L}\Psi_2
\overline{L}^{s-2}T^{m+1} u\|_{L^2})^2 + \varepsilon \|
\overline{L}\Psi_2 L\overline{L}^{s-1}T^mu\|_{L^2}^2$$
for any $\varepsilon >0$ (plus terms with $\Psi_2$
differentiated and corresponding gains.) 

While the last term on the right has one $L$ derivative,
normally not well estimated, we have seen that we are able to
iterate the {\it a priori} estimate as long as one $\overline{L}$
remains. 

To see that this is the case, we begin the iteration of the {\it a
priori} estimate with the function
$v=\Psi L\overline{L}^{s-1}T^mu$ instead of $\Psi
\overline{L}^{s}T^mu.$ The corresponding bracket to
consider is (modulo $\Psi^\prime$):
$$([L \overline{L}, \Psi_2 L\overline{L}^{s-1}]T^m u,
\Psi_2
L\overline{L}^{s-1}T^m u)_{L^2}\sim$$
$$\sim s(\Psi
L\overline{L}^{s-1}T^{m+1} u, \Psi
L\overline{L}^{s-1}T^m u)_{L^2}$$
$$\sim s(\Psi
\overline{L}^{s-1}T^{m+1} u, \overline{L}L\Psi
\overline{L}^{s-1}T^m u)_{L^2}$$
which, after a weighted Schwarz inequality, reduces to the
previous case (iterated one step) and a small multiple of the
left hand side.

To conclude the proof we observe that so far we have passed
from estimating $p$ derivatives in $\Omega_0$ to estimating
$p/2$ derivatives in $\Omega_1:$
$${\|D^pu\|_{L^2(\Omega_0)} \over p!}
\leq (C_0/d_{0})^p\sup_{\tilde{p}\leq
p/2}{\|D^{\tilde{p}}u\|_{L^2(\Omega_1)}
\over \tilde{p}!}$$ 
with $d_{0}$ equal to the separation between $\Omega_0$
and the complement of $\Omega_1.$ 

Lastly we need to nest $\log_2 p$ open sets between
$\Omega_0$ and the set $\tilde\Omega$ where $f$ is assumed
to be analytic with separations $d_{0}, d_{1},
\ldots d_{\log_2(p)}$ in such a way that the sum is bounded
independently of $p$ and the product is less than a universal
constant raised to the power $p,$ e.g., such that 
$$\sum_0^{\log_2(p)} d_j \leq 1,\qquad \prod_0^{\log_2(p)}
d_j^{-p/2^j} \leq C^p.$$

But $d_j$ = $1/2^{j+1}$ will satisfy these conditions.

\end{document}